\documentclass[12pt]{article}
\usepackage{amsfonts}
\usepackage{amsmath}
\usepackage{latexsym}
\usepackage{amscd}

\newtheorem{lemma}{Lemma}

\def\X{{\cal X}}

\newcommand{\be}{\begin{equation}}
\newcommand{\ee}{\end{equation}}

\begin{document}
\bibliographystyle{plain}

\thispagestyle{empty}
\setcounter{page}{0}

\vspace {2cm}

{\Large  Guszt\'av Morvai and Benjamin Weiss : }

\vspace {1cm}

{\Large  Forward Estimation for   Ergodic Time Series }

\vspace {1cm}

{\large Ann. Inst. H. Poincar\'e Probab. Statist.  41  (2005),  no. 5, 859--870.}

\vspace {1cm}

\begin{abstract}
The forward estimation  problem for stationary and ergodic  time series $\{X_n\}_{n=0}^{\infty}$ taking values from a finite alphabet 
${\cal X}$ is to 
estimate the probability that $X_{n+1}=x$ based on the observations 
$X_i$, $0\le i\le n$ without prior knowledge of the distribution of the process $\{X_n\}$. 
 We present a simple procedure $g_n$ which is evaluated on the data segment $(X_0,\dots,X_n)$ and for which, 
 ${\rm error}(n)=|g_{n}(x)-P(X_{n+1}=x|X_0,\dots,X_n)|\to 0$ almost surely for a subclass of all 
 stationary and ergodic time series, while for the full class the Cesaro average of the error tends to zero 
 almost surely 
 and moreover, 
 the error tends to zero in probability. 

Le probl\`eme d'estimation future d'une s\'erie de temps ergodique et
stationnaire $\{X_n\}_{n=0}^{\infty}$, qui prend ses valeures dans un
alphabet fini ${\cal X}$, est d'estimer la probabilit\'e que 
$X_{n+1} = x$, connaissant les $X_i$ pour $0\leq i \leq n$ mais sans
connaissance pr\'ealable de la distribution du processus $\{X_i\}$. 
 Nous pr\'esentons un proc\'ed\'e simple $g_n$, evalu\'e dur les donn\'ees
$(X_0,\dots,X_n)$, pour lequel 
${\rm erreur}(n) = | g_n(x) -P(X_{n+1} = x|X_0,\dots, X_n) \to 0$
presque s\^urement pour une sous-classe de toutes les s\'eries de temps
ergodiques et stationnaires, tandis que pour la classe enti\`ere la
moyenne de Cesaro de l'erreur tend vers z\'ero presque s\^urement. De
plus, l'erreur tend vers z\'ero en probabilit\'e.
\end{abstract}



\pagebreak

\section{Introduction}

T. Cover \cite{Cover75} posed two fundamental problems concerning 
estimation for stationary and ergodic binary  time series $\{ X_n\}_{n=-\infty}^{\infty}$.
(Note that a stationary time series $\{ X_n\}_{n=0}^{\infty}$ 
can be extended to be a two sided stationary time series 
$\{ X_n\}_{n=-\infty}^{\infty}$.) 
Cover's first problem was on  backward estimation. 

\smallskip
\noindent
{\bf Problem 1} {\it Is there an estimation scheme $f_{n}$ for the value 
$P(X_1=1|X_{-n},\dots,X_{0})$ 
such that $f_{n}$ depends solely on the observed data 
segment $(X_{-n},\dots,X_{0})$ 
and 
$$
\lim_{n\to\infty} |f_{n}(X_{-n},\dots,X_{0})-P(X_{1}=1|X_{-n},\dots,X_{0})|=0
$$
almost surely for all stationary and ergodic binary   
time series $\{ X_n\}_{n=-\infty}^{\infty}$? 
}

\smallskip
\noindent
This problem was solved by Ornstein \cite{Ornstein78} by constructing such a scheme. (See also Bailey \cite{Bailey76}.)
Ornstein's scheme is not 
a simple one and  the proof of consistency is rather sophisticated. 
For an even more general case, a much simpler scheme and proof of consistency were provided by Morvai, Yakowitz, Gy\"orfi \cite{MoYaGy96}. 
(See also Algoet \cite{Algoet92} and  Weiss \cite{Weiss00}.) Note that none of thsese  schemes are reasonable from the 
data consumption point of view. 

\smallskip
\noindent
Cover's second problem was on forward estimation.  

\smallskip
\noindent
{\bf Problem 2} {\it Is there an estimation scheme $f_{n}$ for the value 
$P(X_{n+1}=1|X_0,\dots,X_n)$ such that 
$f_{n}$ depends solely on the data segment $(X_0,\dots, X_n)$ and 
$$
\lim_{n\to\infty} |f_{n}(X_0,\dots,X_n)-P(X_{n+1}=1|X_0,\dots,X_n)|=0
$$
almost surely for all stationary and ergodic binary    
time series $\{ X_n\}_{n=-\infty}^{\infty}$? }

\smallskip
\noindent
This problem was answered by Bailey \cite{Bailey76} in a negative way, that is, he showed that there is no such scheme. 
(Also see Ryabko \cite{Ryabko88}, Gy\"orfi,  Morvai, Yakowitz \cite{GYMY98}
 and  Weiss \cite{Weiss00}.) 
Bailey used the technique of cutting and stacking developed by Ornstein \cite{Ornstein74} and 
Shields \cite{Shields91}. Ryabko's construction was  based on a function of 
an infinite state Markov-chain. 

\noindent
Morvai \cite{Mo00} addressed a modified version of Problem 2. There one is not required to predict for all time 
instances
 rather he may refuse to predict for certain values of $n$. 
However, he is expected to predict infinitely often. 
Morvai \cite{Mo00} proposed a sequence of stopping times $\lambda_n$ and he managed to estimate the conditional 
probability $P(X_{\lambda_n+1}=1|X_0,\dots,X_{\lambda_n})$ in the pointwise sense, that is, for his estimator  
 along the proposed stopping time sequence, the error tends to zero as $n$ increases, almost surely. 
 Another estimator  was proposed for this modified Problem 2 by Morvai and Weiss \cite{MW03} for which 
 the $\lambda_n$ grow more slowly, but the consistency only holds   for 
 a certain subclass of all stationary binary time series.

\noindent
In this paper we consider the original Problem 2 but we shall  impose an additional restriction on the possible 
time series.
The conditional probability $P(X_{1}=1|\dots,X_{-1},X_0)$ is said to be 
continuous if a version of it   is continuous   
with respect to metric $\sum_{i=0}^{\infty} 2^{-i-1} |x_{-i}- y_{-i}|$,
 where $x_{-i},y_{-i}\in \{0,1\}$. 

\smallskip
\noindent
{\bf Problem 3} {\it Is there an estimation scheme $f_{n}$ for the value 
$P(X_{n+1}=1|X_0,\dots,X_n)$ such that 
$f_{n}$ depends solely on the data segment $(X_0,\dots, X_n)$ and 
$$
\lim_{n\to\infty} |f_{n}(X_0,\dots,X_n)-P(X_{n+1}=1|X_0,\dots,X_n)|=0
$$
almost surely for all stationary and ergodic binary    
time series $\{ X_n\}_{n=-\infty}^{\infty}$ with continuous 
conditional probability $P(X_{1}=1|\dots,X_{-1},X_0)$? }

\smallskip
\noindent
We will answe this question in the affirmative. 
This class  includes all $k$-step Markov chains.   It is not known if the schemes proposed 
by  Bailey \cite{Bailey76}, Ornstein \cite{Ornstein78}, 
 Morvai, Yakowitz, Gy\"orfi \cite{MoYaGy96} solve Problem 3 or not.

\smallskip
\noindent
{\bf Problem 4}
{\it Is there an estimation scheme $f_{n}$ for the value 
$P(X_{n+1}=1|X_0,\dots,X_n)$ such that 
$f_{n}$ depends solely on the data segment $(X_0,\dots, X_n)$ and 
$$
\lim_{n\to\infty} {1\over n} \sum_{i=0}^{n-1} |f_i(X_0,\dots,X_i)- P(X_{i+1}=1|X_0,\dots,X_i)|=0 
$$
almost surely for all stationary and ergodic binary    
time series $\{ X_n\}_{n=-\infty}^{\infty}$? }

\smallskip
\noindent
Bailey \cite{Bailey76} (cf.  Algoet \cite{Algoet94} also) showed that any scheme that 
solves Problem 1 can 
be easily modified to solve Problem 4 (indeed, just exchange the data segment $(X_{-n},\dots,X_0)$ for 
$(X_0,\dots,X_n)$, but apparently not all solutions of Problem 4 arise in this fashion. 
For further reading cf. Algoet \cite{Algoet92}, \cite{Algoet99}, Morvai, Yakowitz, Gy\"orfi \cite{MoYaGy96},
 Gy\"orfi et. al. \cite{GYKKW02},  Gy\"orfi, Lugosi and Morvai \cite{GYLM99},
   Gy\"orfi and Lugosi \cite{GYL02} and  
 Weiss \cite{Weiss00}.

\smallskip 
\noindent
{\bf Problem 5}
{\it Is there an estimation scheme $f_{n}$ for the value 
$P(X_{n+1}=1|X_0,\dots,X_n)$ such that 
$f_{n}$ depends solely on the data segment $(X_0,\dots, X_n)$ and 
 for arbitrary $\epsilon>0$,
$$
\lim_{n\to\infty} P( |f_i(X_0,\dots,X_i)- P(X_{i+1}=1|X_0, \dots,X_i)|>\epsilon)=0.
$$
for all stationary and ergodic binary    
time series $\{ X_n\}_{n=-\infty}^{\infty}$? }

\smallskip
\noindent
By stationarity, for any scheme that solves Problem 1, the shifted version of it solves Problem 5.
(Just replace the data segment $(X_{-n},\dots,X_0)$ by $(X_0,\dots,X_n)$.)

\bigskip
\noindent
There are existing schemes that solve Problem 4 (e.g. Bailey \cite{Bailey76}, Ornstein \cite{Ornstein78}, 
and even for a more general case Morvai, Yakowitz, Gy\"orfi \cite{MoYaGy96}, 
Algoet \cite{Algoet92}, Gy\"orfi and Lugosi \cite{GYL02})
 and there are schemes that solve Probkem 5 (e.g.  Bailey \cite{Bailey76}, Ornstein \cite{Ornstein78} and for even more general case 
  Morvai, Yakowitz, Gy\"orfi \cite{MoYaGy96}, 
Algoet \cite{Algoet92}, Morvai, Yakowitz and Algoet \cite{MoYaAl97}). 
In this paper we  propose a reasonable, very simple  algorithm that {\it simultanously} solves Problem 3, 4 and 5.  
 Note that the schemes given by  Bailey \cite{Bailey76}, Ornstein \cite{Ornstein78}, 
 Morvai, Yakowitz, Gy\"orfi \cite{MoYaGy96},  Algoet \cite{Algoet92} and Weiss \cite{Weiss00}
are not reasonable at all, they  consume data extremely  rapidly, cf. Morvai \cite{MorvaiPhD} and it is 
not known if their schemes solve Problem ~3 or not.

\section{Preliminaries and Main Results}

Let $\{X_n\}_{n=-\infty}^{\infty}$ be a stationary  time series taking values from a 
finite alphabet 
${\cal X}$. (Note that all stationary time series $\{X_n\}_{n=0}^{\infty}$ 
can be thought to be a 
two sided time series, that is, $\{X_n\}_{n=-\infty}^{\infty}$. )  
For notational convenience, let $X_m^n=(X_m,\dots,X_n)$,
where $m\le n$. Note that if $m>n$ then $X_m^n$ is the empty string. 

\noindent
Let $g: {\cal X}\rightarrow (-\infty,\infty)$ be  arbitrary. 

\noindent
Our goal is to estimate the conditional expectation $E(g(X_{n+1})|X_0^n)$ from samples 
$X_0^n$.

\noindent
For $k\ge 1$  define the stopping times  
$\tau^{k}_i(n)$ which indicate where the $k$-block $X_{n-k+1}^{n}$ occurs previously in the time series $\{X_n\}$. 
Formally we set $\tau^{k}_0(n)=0$ and for $i\ge 1$ let 
\begin{equation}\label{deftau} 
\tau^k_i(n)= \min\{t>\tau_{i-1}^k(n) : X_{n-k+1-t}^{n-t}=X_{n-k+1}^{n}\}.  
\end{equation}

\noindent
Let $K_n\ge 1$ and $J_n\ge 1$ be sequences of  nondecreasing positive integers tending to $\infty$ which will be fixed later.  

\noindent
 Define $\kappa_n$ as the largest $1\le k\le K_n$ such that there are at least 
$J_n$ occurrences of  the block 
$X^n_{n-k+1}$ in the data segment $X_0^n$, that is, 
\begin{equation}\label{defkappa}
 \kappa_n=\max\{ 1\le k\le \ K_n : \tau^{k}_{J_n}(n) \le n-k+1\} 
\end{equation}
 if
there is such $k$ and $0$ otherwise. 

\noindent
Define $\lambda_n$ as the number of occurences of the block $X^n_{n-\kappa_n+1}$ 
in the data segment $X_0^n$,  that is,  
\begin{equation}\label{deflambda}
\lambda_n=\max\{1\le j: \ \tau_j^{\kappa_n}\le n-\kappa_n+1\}
\end{equation}
if $\kappa_n>0$  and zero otherwise. 
Observe that if $\kappa_n>0$ then $\lambda_n\ge J_n$. 

\noindent
Our estimate $g_n$ for $E(g(X_{n+1})|X_0^n)$ is defined as $g_0=0$ and for $n\ge 1$, 
\begin{equation}\label{estimatordef} 
g_n={1\over \lambda_n }\sum_{i=1}^{\lambda_n} g(X_{n-\tau^{\kappa_n}_i(n)+1})
\end{equation}
if $\kappa_n>0$ and zero otherwise. 

\bigskip
\noindent
Let ${\cal X}^{*-}$ be the set of all one-sided   sequences, that is, 
$${\cal X}^{*-} =\{ (\dots,x_{-1},x_0): x_i\in {\cal X} \ \ \mbox{for all $-\infty<i\le 0$}\}.$$
Define the function 
$G : {\cal X}^{*-}\rightarrow (-\infty,\infty)$ 
as 
$$
G(x^{0}_{-\infty})=E(g(X_1)|X^{0}_{-\infty}=x^0_{-\infty}).
$$
Note that as a conditional expectation this is only defined almost surely.  
E.g. if $g(x)=1_{\{x=z\}}$ for a fixed 
$z\in {\cal X}$ then $G(y^0_{-\infty})=P(X_1=z|X_{-\infty}^0=y^0_{-\infty}).$  

\noindent
Define a distance on ${\cal X}^{*-}$ as   
$$
d^*(x^0_{-\infty},y^0_{-\infty})=
\sum_{i=0}^{\infty} 2^{-i-1} 1_{ \{x_{-i}\neq y_{-i}\} }.
$$

\smallskip
\noindent
{\sc Definition}
The conditional expectation  $G(X^{0}_{-\infty})$
is said to be continuous if a version of it   is continuous on the set ${\cal X}^{*-}$  
with respect to metric $ d^*(\cdot,\cdot)$. 
Since this space is compact, in fact, continuity is equivalent to uniform continuity. 

\bigskip
\noindent
The processes with continuous conditional expectation are essentially the Random Markov Processes of Kalikow~ 
\cite{Ka90}, or the continuous g-measures studied by Mike Keane  \cite{Ke72}.

\bigskip
\noindent 
{\bf Theorem}
{\it
Let $\{X_n\}$ be a stationary  and ergodic 
time series taking values from a finite alphabet $\cal X$. 
Assume  $K_n=\max(1, \lfloor 0.1 \log_{|{\cal X}|} n \rfloor) $ and
 $J_n=\max(1,\lceil n^{0.5}\rceil)$. Then

\noindent
{\bf (A)} if the conditional expectation  
$G(X_{-\infty}^{0})$ 
is  continuous with respect to metric $d^*(\cdot,\cdot)$ then  
$$
\lim_{n\to\infty} \left| g_n- E(g(X_{n+1})|X_0^{n}) \right| =0\ \ \mbox{almost surely,}
$$
 
\noindent
{\bf (B)} without any continuity assumption, 
$$
\lim_{n\to\infty} {1\over n} \sum_{i=0}^{n-1} |g_i- E(g(X_{i+1})|X_0^{i})|=0\ \ \mbox{almost surely,} 
$$
{\bf (C)} without any continuity assumption, for arbitrary $\epsilon>0$,
$$
\lim_{n\to\infty} P( |g_n- E(g(X_{n+1})|X_0^{n})|>\epsilon)=0.
$$
}

\bigskip
\noindent
{\bf Remarks:}

\noindent
 Note that these results are valid vithout the ergodic assumption. 
One may use the ergodic decomposition throughout the proofs,  
cf. Gray \cite{Gr88} p. 268.

\bigskip
\noindent
We note that  from the proof of Ryabko~\cite{Ryabko88} and Gy\"orfi, Morvai, Yakowitz~\cite{GYMY98} it is clear 
that the continuity condition in the first part of the Theorem can not be relaxed. Even for  
the class of all stationary and ergodic  binary 
time-series with merely almost surely continuous conditional probability 
$P(X_1=1|\dots,X_{-1},X_{0})$ one can not solve Problem 2 in the Introduction. 
(An almost surely continuous  conditional probability is such that as a function restricted to a set $C$ with 
full measure, it is continuous
on $C$. )

\bigskip
\noindent
We do  not know if the shifted version of our proposed scheme $g_n$ solves Problem~1 or not.  
(That is, in the case when $g_n$ is evaluated on $(X_{-n},\dots,X_0)$ rather than on $(X_0,\dots,X_n)$. 
  
\bigskip
\noindent
If $\cal X$ is a countably infinite alphabet then there is no scheme that could achieve similar 
result to part (A) in the Theorem for all bounded $g(\cdot)$, 
even if you assume that the resulting $G(\cdot)$ is continuous, and the time series is in fact 
a first order Markov chain.  
Indeed, whenever a new state appears  which has not occured before, you are unable to predict, 
cf.  Gy\"orfi,  Morvai, Yakowitz \cite{GYMY98}.  

\section{Auxiliary Results} 

\noindent
For $k\ge 1$, $n\ge 0$ and $j\ge 0$ it will be useful to define auxiliary  processes  
$\{ {\tilde X}^{(k,n,j)}_i\}_{i=-\infty}^{\infty}$ as follows.
Let 
\begin{equation} \label{deftildex}
{\tilde X}^{(k,n,j)}_i=X_{n-\tau^k_{j}(n)+i} \ \ \mbox{for $-\infty< i<\infty$.}
\end{equation}
For an arbitrary stationary time series $\{Y_n\}$ for $k\ge 1$ let 
 $\tilde\tau^k_0(Y_{-\infty}^{\infty})=0$ and for $i\ge 1$ define   
\begin{equation}\label{deftildetau}
{\tilde\tau}^k_i(Y_{-\infty}^{\infty})=
\min\{t>{\tilde \tau}^k_{i-1}(Y_{-\infty}^{\infty}) : Y_{-k+1+t}^{t}=Y_{-k+1}^{0}\}.
\end{equation}
If it is obvious on which time series ${\tilde\tau}^k_i(Y_{-\infty}^{\infty})$ is evaluated, we will write ${\tilde\tau}^k_i$.

\noindent
Let $T$ denote the left shift, that is, $(T x^{\infty}_{-\infty})_i=x_{i+1}$. 

\bigskip
\noindent
We will need the next lemmas for later use.

\begin{lemma} \label{equaldistr1}
Let $\{X_n\}$ be a stationary time series taking values from a finite alphabet $\cal X$. 
For $k\ge 1$, $n\ge 0$ and  $j\ge 0$,  the time series 
$\{ {\tilde X}^{(k,n,j)}_{i}\}_{i=-\infty}^{\infty}$ 
has the same distribution as $\{ X_i\}_{i=-\infty}^{\infty}$. 
\end{lemma}
{\sc Proof}
Note that by (\ref{deftau}),  and (\ref{deftildetau}), 
$$
T^{n-s} 
\{ X_{n-\tau^k_{j}(n)+l}^{n-\tau^k_{j}(n)+m}=x_l^m,
\tau^k_{j}(n)=s\}=
\{ X_l^m=x_l^m,{\tilde\tau}^k_{j}=s\}
$$
where ${\tilde \tau}^k_j$ is evaluated on time series $\{X_i\}_{i=-\infty}^{\infty}$.
Now by (\ref{deftildex}) and stationarity,  
\begin{eqnarray*}
 P({\tilde X}^{(k,n,j)}_l=x_l, \dots, {\tilde X}^{(k,n,j)}_m=x_m)
&=& \sum_{s=0}^{\infty}P( X_{n-\tau^k_{j}(n)+l}^{n-\tau^k_{j}(n)+m}=x_l^m,
\tau^k_{j}(n)=s)\\
&=&\sum_{s=0}^{\infty}P( X_l^m=x_l^m,{\tilde\tau}^k_{j}=s)\\
&=&P(X_l^m=x_l^m)
\end{eqnarray*}
and the proof of Lemma~\ref{equaldistr1} is complete.

\smallskip
\begin{lemma} \label{kntoinftylemma}
Let $\{X_n\}$ be a stationary and ergodic time series taking values from a finite alphabet $\cal X$. 
Assume $K_n \to\infty$, $J_n \to\infty$ and 
$\lim_{n\to\infty} {J_n\over n}=0$. 
Then  
$$
\lim_{n\to\infty} \kappa_n=\infty \ \ \mbox{almost surely.}
$$

\end{lemma}
{\sc Proof}
We argue by contradiction. 
Suppose, that $\kappa_{n_i}=K$, $X^{n_i}_{n_i-K}=x^0_{-K}$ for a subsequence $n_i$.
Then a simple frequency count
(in the data segment $X_0^{n_i}$ there are less than   $J_{n_i}$ 
occurrences of block $x^0_{-K}$) yields
that   
$$P(X^0_{-K}=x^0_{-K})\le \lim_{n\to\infty}  { J_n\over n}= 0.$$
The set of sequences that contain a block with zero probability has zero probability and thus   
 Lemma~\ref{kntoinftylemma} is proved.

\section{Pointwise Consistency}

\smallskip

{\sc Proof of  Theorem (A).}
By Lemma~\ref{kntoinftylemma}, for large n, 
\begin{eqnarray*}
\lefteqn{|g_n(x) -E(g(X_{n+1})|X_0^{n})|}\\
&=&\left| {1\over \lambda_n }\sum_{j=1}^{\lambda_n } 
g(X_{n-\tau^{\kappa_n}_j(n)+1})  -E(g(X_{n+1})|X_0^{n})   \right| \\
&\le& \max_{J=J_n,\dots,n} \max_{k=1,\dots,K_n }
\left| {1\over J }\sum_{j=1}^{J } [g( X_{n-\tau^k_j(n)+1})    -
G(X_{-\infty}^{n-\tau^k_j(n)}) ] \right| \\
&+& \left| {1\over \lambda_n } \sum_{j=1}^{\lambda_n }
 G(X_{-\infty}^{n-\tau^{\kappa_n}_j(n)})
 -E(g(X_{n+1})|X_0^{n}) \right| .
\end{eqnarray*}
Concerning the first term, by (\ref{deftau}), ({\ref{deftildetau}) and  (\ref{deftildex}), 
\begin{eqnarray} \nonumber
{1\over J }\sum_{j=1}^{J }
\lefteqn{ [g(X_{n-\tau^k_j(n)+1}) -G(X_{-\infty}^{n-\tau^k_j(n)}) ]}\\
&=& 
\label{aux1}
{1\over J }\sum_{j=0}^{J -1}
 [g(\tilde X^{(k,n,J)}_{\tilde\tau^k_j+1}) -
G(\dots,\tilde X^{(k,n,J)}_{\tilde\tau^k_j-1},\tilde
X^{(k,n,J)}_{\tilde\tau^k_j})]
\end{eqnarray} 
where $\tilde\tau^k_j$ is evaluated on $\{\tilde X^{(k,n,J)}_i\}_{i=-\infty}^{\infty}$.
Since by Lemma~\ref{equaldistr1} 
$$
G(\dots,\tilde X^{(k,n,J)}_{\tilde\tau^k_j-1},\tilde X^{(k,n,J)}_{\tilde\tau^k_j})=
E(g({\tilde X}^{(k,n,J)}_{\tilde\tau^k_j+1})|
\dots,\tilde X^{(k,n,J)}_{\tilde\tau^k_j-1},\tilde X^{(k,n,J)}_{\tilde\tau^k_j}), 
$$ 
the pair 
 $\left(\Gamma_j= g(\tilde X^{(k,n,J)}_{\tilde\tau^k_j+1})
-G(\dots,\tilde X^{(k,n,J)}_{\tilde\tau^k_j-1},\tilde X^{(k,n,J)}_{\tilde\tau^k_j}), 
{\cal F}_j=\sigma(X_{-\infty}^{\tilde\tau^k_j})\right)
$ 
 forms a martingale difference sequence 
( $E(\Gamma_j|{\cal F}_j)=0$ and $\Gamma_j$ is measurable with respect to ${\cal F}_{j+1}$) 
for which Azuma's exponential bound (cf. Azuma~\cite{Azuma67})
yields
$$
P\left(\left| {1\over J }\sum_{j=0}^{J-1 }
 [g(\tilde X^{(k,n,J)}_{\tilde\tau^k_j+1})-G(\dots,\tilde X^{(k,n,J)}_{\tilde\tau^k_j})]
\right|>\epsilon\right)\le 2 e^{-\epsilon^2 J \over B}
$$
for any $B$ such that  $\max_{x\in \X} |g(x)|<B$. 
Now by (\ref{aux1})
$$
P\left(\max_{J=J_n,\dots,n} \max_{k=1,\dots,K_n}
\left|{1\over J } \sum_{j=1}^{J }
 [ g( X_{n-\tau^k_j(n)+1})-G(X_{-\infty}^{n-\tau^k_j(n)})
\right|>\epsilon\right)\le n K_n  2 e^{-\epsilon^2 J_n \over B}
$$
and by assumption 
$n K_n  2 e^{-\epsilon^2 J_n/B}$
 sums up  and the Borel-Cantelli Lemma yields almost sure convergence to zero. 
Concerning the second term,  
$$\left|{1\over  \lambda_n } \sum_{j=1}^{\lambda_n } 
 [G(X_{-\infty}^{n-\tau^{\kappa_n}_j(n)})
 -E\left( G(X_{-\infty}^{n})|X^n_{0} \right)  ]\right|\to 0 \ \ \mbox{almost surely}
$$
since $\kappa_n$ tends to infinity by Lemma~\ref{kntoinftylemma}, 
$X^{n-\tau_j^{\kappa_n}(n)}_{n-\tau_j^{\kappa_n}(n)-\kappa_n+1}=X^n_{n-\kappa_n+1}$ for 
$0\le j\le \lambda_n$, 
and the conditional expectation $G(\cdot)$  is  in fact uniformly continuous on ${\cal X}^{*-}$ 
with respect to 
$d^*(\cdot,\cdot)$.  The
proof
of Theorem (A)  is complete.

\section{Time Average Performance}

\smallskip
\noindent
If the process does not have continuous conditional expectations then the last step in the proof of 
Theorem (A) is not valid. It can be carried out for most time instances $n$ by using the typical  behaviour of almost every realization 
$x_{-\infty}^{\infty}$. 
More specifically, for every  $\delta>0$, the probability of the  set of those $ x_{-\infty}^{0}$ for which 
$$|E(g(X_1|X^0_{-k+1}=X^0_{-k+1})-G(X^0_{-\infty})|<\delta \ \ \mbox{for all $k\ge K$}$$  
tends to one as $K$ tends to infinity. 
The typical behaviour we are after is the statement that most of  the times  $t=n-\tau_j^{\kappa_n}(n)$ 
the sequence $T^t x^t_{-\infty}$ 
belongs to the above mentioned  set. While this need not be the case for all $n$, it is true for most $n$'s and 
the next lemma makes this precise. For the analysis we will fix a value of $\kappa_n$ at $k$. 
 
 \noindent
 Define the set of good indexes $M_n(\delta,K)\subseteq \{ K-1,\dots, n-1\} $ as 
\begin{eqnarray*}
\lefteqn{M_n(\delta,K)}\\
&=&\{K-1\le i\le n-1 \ :  \ |E(g(X_{i+1})|X^i_{i-k+1})-G(X^i_{-\infty})|
<\delta \ \ \mbox{for all $k\ge K$} \}.
\end{eqnarray*}

\noindent
We will analyze the behaviour of our algorithm for $\kappa_n=k$ for each $i\le n$ by first dividing up 
the indices $\{1,2,\dots, n\}$ according to the value of $X^i_{i-k+1}=y^0_{-k+1}$, and considering what happens 
for each of these. 

\noindent
Let $y^0_{-k+1}\in {\cal X}^k$. Define the set of indexes $I_n^k(y^0_{-k+1})\subseteq 
\{ k-1,\dots, n-1\} $, where you can find the pattern
$y^0_{-k+1}$, that is,  
$$
I_n^k(y^0_{-k+1})=\{k-1\le i\le n-1 \ : \  X^i_{i-k+1}=y^0_{-k+1}\}.
$$
Define $D_k(i)$ as
\begin{equation*}
D_k(i)=\left\{
\begin{array}{ll}
\{\tau^k_j(i)   :   \mbox{$\tau^k_j(i)\le i-k+1$ and $1 \leq j\leq i+1$\}} \   \ &  \mbox{if $\tau^k_{J_i}(i)\le i-k+1$} \\
\emptyset  &  \mbox{otherwise.} 
\end{array}
\right.
\end{equation*}
Let  $E_n^k(\delta,K)$ be defined as
$$
E_n^k(\delta,K)=\{0\le i\le n-1: \ \ 
\mbox{$|D_k(i)\bigcap M_n(\delta,K)|> (1-{\delta}^{0.5})|D_k(i)|$.}\}
$$

\noindent
If the number of occurrences of $y^0_{-k+1}$ prior to $i$ was not enough for our algorithm then $D_k(i)$ 
will be empty. This is rare, and can be expressed as follows: Let 
$$
F_n^k=\{0\le i\le n-1: \ \ \mbox{$D_k(i)=\emptyset$.}\}
$$
It is immediate that 
\begin{equation} \label{aux2}
|F_n^k|\le |{\cal X}|^k J_n. 
\end{equation}

\begin{lemma} \label{mostofthetimeinElemma} 
Assume $|M_n(\delta,K)|\ge (1-\delta)n$.  
Then 
$$
\left| \{ 0\le l\le n-1 \ : \ l\notin E_n^k(\delta,K) \ \mbox{and} \ l\notin F_n^k\}
\right| 
\le {\delta}^{0.5} n.
$$
\end{lemma} 
{\sc Proof} 
Fix $\delta$, $K$, $k$ and $x\in {\cal X}$. Temporarily
fix also $y^0_{-k+1}\in {\cal X}^k$. 
Let $z=|I_n^k(y^0_{-k+1})| $ and let 
$k\le i_1\le i_2\le \dots\le i_z$ 
denote the elements of
$I_n^k(y^0_{-k+1})$. 
Let $i_j(y^0_{-k+1})$ be the largest element $i_{j\prime}$ of $I_n^k(y^0_{-k+1})$ such that 
$D_k(i_{j\prime})\neq \emptyset$ and 
$$|\{ 0\le l\le n-1 \ : \ 
l\in D_k(i_{j\prime}) \ \mbox{and} \ l\notin M_n(\delta,K) \}|\ge
{\delta}^{0.5} |D_k(i_{j\prime})|.
$$
Define $S$ to be the set of these indexes as $y^0_{-k+1}$  varies over all element ${\cal X}^k$. 
It is clear that if $i,j\in S$, $i\neq j$ then $D_k(i)\bigcap D_k(j)=\emptyset$ since different blocks $y^0_{-k+1}$ are involved. 
It follows from the construction  that  
$\{ D_k(i)\}_{i\in S}$ is a disjoint cover
of 
$\{ 0\le l\le n-1 \ : \ l\notin E_n^k(\delta,K) \ \mbox{and} \ l\notin F_n^k\}$. It follows that 
\begin{eqnarray*} 
n \delta &\ge& 
\left|  \{ 0\le l\le n-1\ : \ l\notin  M_n(\delta,K)\} \right|\\
&\ge& 
\sum_{i\in S} \left| 
\{ 0\le l\le n-1 \ : \ 
l\in D_k(i) \ \mbox{and} \ l\notin M_n(\delta,K) \}
\right| \\ 
&\ge& 
{\delta}^{0.5} \sum_{i\in S}  |D_k(i)| = {\delta}^{0.5} 
|\bigcup_{i\in S} D_k(i)|. 
\end{eqnarray*}
Now 
$$
\left| \{ 0\le l\le n-1 \ : \ l\notin E_n^k(\delta,K) \ \mbox{and} \ l\notin F_n^k\}
\right|\le |\bigcup_{i\in S} D_k(i)| 
\le {\delta}^{0.5} n
$$
and the proof of Lemma~\ref{mostofthetimeinElemma} is complete.

\bigskip
\noindent
{\sc Proof of Theorem (B).} 
Consider
\begin{eqnarray*} 
\lefteqn{
 {1\over n} \sum_{i=0}^{n-1} |g_i(x)- E(g(X_{i+1})|X_0^{i})|}\\
&\le& 
 { |0-E(g(X_{1})|X_0)| \over n}\\
&+&  
 {1\over n} \sum_{i=1}^{n-1}
1_{\{\kappa_i<K_i \}}\\
&+&  {1\over n} \sum_{i=1}^{n-1}
\max_{J=J_i,\dots,i} \left| {1\over J }\sum_{j=1}^{J }
 [g(X_{i-\tau^{K_i}_j(i)+1})   -G(X_{-\infty}^{i-\tau^{K_i}_j(i)}) ] \right|\\
&+& 
 {1\over n} \sum_{i=1}^{n-1}
\left| {1\over  \lambda_i } \sum_{j=1}^{\lambda_i }
 G(X_{-\infty}^{i-\tau^{K_i}_j(i)})-
E(g(X_{i+1})|X^i_{i-K_i+1}) \right|1_{\{\kappa_i=K_i\}} \\
&+&
 {1\over n} \sum_{i=1}^{n-1}
\left|E(g(X_{i+1})|X^i_{i-K_i+1}) 
 -E(g(X_{i+1})|X^i_{0}) \right|.
\end{eqnarray*}
The first term tends to zero. 
The second term tends to zero   
since by (\ref{aux2}) $|F_n^{K_n}|/n\le |{\cal X}|^{K_n}J_n/n \to 0$.

\noindent
Concerning the third  term, by (\ref{aux1})
and by Azuma's exponential bound (cf. Azuma~\cite{Azuma67})  
$$
\sum_{J=J_i}^i P\left(\left| {1\over J }\sum_{j=1}^{J }
 [g({\tilde X}^{(K_i,i,J)}_{\tilde\tau^{K_i}_j+1})  -
G(\dots,{\tilde X}^{(K_i,i,J)}_{\tilde\tau^{K_i}_j-1},
{\tilde X}^{(K_i,i,J)}_{\tilde\tau^{K_i}_j})\} \right| >\epsilon\right)
\le
2i e^{-\epsilon^2 J_i \over B}  
$$ 
(where $B$ is any real such that $2\max_{x\in {\cal X}} |g(x)|<B$)  
and the right hand side is summable, hence  the Borel-Cantelli Lemma yields almost sure convergence to zero.
By Toeplitz lemma the average also converges to zero.  

\noindent
Now we deal with the fourth  term. Let $0< \epsilon<1 $ be arbitrary. 
Choose the integer $d$
large enough such that $|{\cal X}|^{-10 (d-1)}<\epsilon$.  Let $\delta={\epsilon\over d^2 }$.  
Let $K$ and  $N_0$ be so large that 
${|M_n(\delta,K)|\over n}>(1-\delta)$ for all $n\ge N_0$. (There exist such $K$ and $N_0$ since by the ergodic theorem and the martingale convergence theorem
$\lim_{k\to\infty}\lim_{n\to\infty} {|M_n(\delta,k)|\over n}=1$ almost surely.) 
Now let $N_1\ge N_0$ be so large that  $K_n-d+2\ge K$ and 
$|{\cal X}|^{10(K_n-d+1)}\ge N_0$
 for all $n\ge N_1$. 
Assume $n\ge N_1$. 
The sum 
$$
 {1\over n}\sum_{i=1}^{n-1}
\left| {1\over  \lambda_i } \sum_{j=1}^{\lambda_i }
 G(X_{-\infty}^{i-\tau^{K_i}_j(i)})-E(g(X_{i+1})|X^i_{i-K_i+1})
\right|1_{\{ \kappa_i=K_i\} } 
$$
that we are trying to estimate will be divided into blocks according to the value of $K_i$. In fact only values 
in the range $[K_n-d+2,K_n]$ need be considered  since the sum up to $|{\cal X}|^{10 K_n-d+1}$ 
can be estimated  by $|{\cal X}|^{10 (K_n-d+1)}2\max_{y\in {\cal X} } |g(y)|  $ and so by our assumption on $d$,
after dividing by $n$ this will be at most $\epsilon 2\max_{y\in {\cal X} } |g(y)| $. 
For $i$ in the range $[|{\cal X}|^{10 (k-1)}, |{\cal X}|^{10 k})$ 
for $K_n-\delta+2\le k\le K_n$, and $\kappa_i=K_i$,
if $i\in  E_{|\X|^{10 k}}^{k-1} (\delta,K)$ then we get for more than $(1-\sqrt\delta)|D_k(i)|$ terms an 
upper bound of $\delta$
while for the rest we may use $2\max_{y\in {\cal X} } |g(y)|$. This gives an upper bound of 
$$
{\delta |D_k(i)|
  +\sqrt{\delta} |D_k(i)| 2 \max_{y\in {\cal X} } |g(y)|\over |D_k(i)|}.
$$
Using Lemma~\ref{mostofthetimeinElemma} we can estimate the sum over all $i$ in the interval 
$[|{\cal X}|^{10 (k-1)}, |{\cal X}|^{10 k})$ by 
$$
n(\delta +\sqrt{\delta} 2 \max_{y\in {\cal X} } |g(y)|)
  + \sqrt\delta n 2 \max_{y\in {\cal X} } |g(y)|.
$$
Dividing by $n$, we have an upper bound:
$$
 \delta +\sqrt{\delta} 2 \max_{y\in {\cal X} } |g(y)|
  + \sqrt\delta  2 \max_{y\in {\cal X} } |g(y)|.
$$
The same argument yields the same upper bound  for the $i$'s in the range $[|{\cal X}|^{10 K_n},n)$.

\noindent
Summing over $k$ in the range $[K_n-d+2,K_n+1]$ yields an upper bound: 
$$
 d\delta +d\sqrt{\delta} 2 \max_{y\in {\cal X} } |g(y)|
  + d\sqrt\delta  2 \max_{y\in {\cal X} } |g(y)|.
$$ 
Recall that $\sqrt\delta d=\sqrt\epsilon$ and this yields an upper bound:
$$
\epsilon +\sqrt\epsilon  2 \max_{y\in {\cal X} } |g(y)|
  + \sqrt\epsilon  2 \max_{y\in {\cal X} } |g(y)|.
$$
Since $\epsilon$ was arbitrary, the fourth term tends to zero. 

\noindent
Now we deal with the last term. Since by the martingale convergence theorem,
$ E(g(X_{1})|X^0_{-i})\to G(X^0_{-\infty})$ almost surely, thus 
$$
\lim_{i\to\infty} |E(g(X_{1})|X^0_{-K_i+1})-E(g(X_{1})|X^0_{-i})|=0
$$
and applying Breiman's generalized ergodic theorem, cf. Maker \cite{Ma40} (or Algoet \cite{Algoet94}), 
$$
\lim_{n\to\infty} {1\over n} \sum_{i=0}^{n-1}
\left|E(g(X_{i+1})|X^i_{i-K_i+1}) 
 -E(g(X_{i+1})|X^i_{0})\right|=0  
$$
almsost surely and 
the proof of Theorem (B) is complete.

\section{Weak Consistency}

\bigskip
\noindent
{\sc Proof of Theorem (C).} 

\bigskip
\noindent
In order to show that for all ergodic stationary processes our estimate $g_n$ converges in probability we follow the 
steps in the proof of Theorem~(A).
The probability that 
$$
\left( |g_n(x) -E(g(X_{n+1})|X_0^{n})|> 3\epsilon\right)
$$
can be estimated as the sum of the probability of several sets,
$$
 P\left( \max_{J=J_n,\dots,n} \max_{k=1,\dots,K_n }
\left| {1\over J }\sum_{j=1}^{J } [g( X_{n-\tau^k_j(n)+1})    -
G(X_{-\infty}^{n-\tau^k_j(n)}) ] \right| > \epsilon\right),
$$

$$
P(\kappa_n< K_n),
$$

$$
P(\left|E(g(X_{n+1})|X_0^{n}) -E(g(X_{n+1})|X_{n-K_n+1}^{n})\right| >\epsilon )
$$
and

$$
P\left( \left| {1\over \lambda_n } \sum_{j=1}^{\lambda_n }
 G(X_{-\infty}^{n-\tau^{\kappa_n}_j(n)})
 -E(g(X_{n+1})|X_{n-\kappa_n+1}^{n}) \right| > \epsilon,\kappa_n=K_n\right).
$$
For the first, the argument given there suffices. Concerning the second, it tends to zero by  
Lemma \ref{applemma}  in the Appendix.(Apply it with $A=\{X^n_{n-K_n+1}= x^n_{n-K_n+1}\}$, $D=J_n$. 
Then sum over 
all possible $ x^n_{n-K_n+1}$ to get that this second  probabilty in question is not greater than 
$|{\cal X}|^{K_n} J_n/n$ which tends to zero.) 
 For the  third, it is easy to see that 
it tends to zero  by stationarity and  by the  martingale convergence theorem which implies that 
$$
\lim_{n\to \infty}P(\left|E(g(X_{1})|X_{-n}^{0}) -E(g(X_{1})|X_{-K_n+1}^{0})\right| >\epsilon )=0.
$$  
We concentrate on the last probability. 
Recall the notations from the proof of Theorem~(B). The main thing is to show that with probability at least $1-\epsilon$, for $n$ sufficiently large, 
most of the elements $l\in  I_n^{K_n}(X^0_{-K_n+1})$ are such that $T^l x^{\infty}_{-\infty}$ does not 
belong to the set 
 $${\tilde M}_n(\epsilon)=
  \{ x^{\infty}_{-\infty} \ :  \ |E(g(X_{1})|X^0_{-k+1}=x^0_{-k+1})-G(x^0_{-\infty})|
>\epsilon \ \ \mbox{for some $k\ge K_n$} \}
$$
as neither does $T^n x^{\infty}_{-\infty}$ itself. By the martingale convergence theorem, the probability of the set 
${\tilde M}_n(\epsilon)$ tends to zero as $n$ tends to infinity. 
Let $n$ be so large that this  probability in question is less than $\epsilon^2/2$. 
Let 
$$
B_n=\{x^{\infty}_{-\infty}:  |\{l\in I_n^{K_n}(x^n_{n-K_n+1}): 
x^{\infty}_{-\infty} \in T^{-l}{\tilde M}_n(\epsilon)\}| 
>\epsilon| I_n^{K_n}(x^n_{n-K_n+1})|\}.
$$
The probability of $B_n$ will be evaluated using the ergodic theorem along the orbit of a typical point.  
Let $x_{-\infty}^{\infty}$ be  such a typical orbit and  $N$ be a very large number.  
Fix $y^0_{-K_n+1}$, and note those elements in $I_N^{K_n}(y^0_{-K_n+1})$  that belong to $B_n$.  
We will cover them with disjoint blocks of length $K_n$, begining on the right end $N-1$ in the obvious way. 
These sets (subsets of  $I_N^k(y^0_{-K_n+1})$ ) we call $C_r(y^0_{-K_n+1})$ where $r=1,2,\dots$. 
Formally, let $\dots<l_2<l_1$ denote the elements of $I_N^{K_n}(y^0_{-K_n+1})$. 
Let $C_0(y^0_{-K_n+1})=\emptyset$.
For $r\ge 1$ we define $C_r(y^0_{-K_n+1})$ recursively. Let $l$ be the largest index such that  
$l\ge n$, $l\not \in \bigcup_{{r\prime}<r}C_{r\prime}(y^0_{-K_n+1})$
and $x^{\infty}_{-\infty}\in T^{-n+l} B_n$. If there is such $l$ then  set 
$C_r(y^0_{-K_n+1})=\{l-n+K_n-1\le l_i\le l\ \ \mbox{for $i=1,2,\dots$} \}$. Let $R(y^0_{-K_n+1})$ be the largest $r$ 
for which $C_r(y^0_{-K_n+1})$ is defined. 
Let $$I_N({\tilde M}_n(\epsilon))=\{ 0\le l\le N-1: T^lx_{-\infty}^{\infty}\in {\tilde M}_n(\epsilon)\}.$$
Then by the construction of $ C_r(y^0_{-K_n+1})$, for each $1\le r\le R(y^0_{-K_n+1})$, 
$$
| \{ l\in C_r(y^0_{-K_n+1}): T^l x_{-\infty}^{\infty} \in {\tilde M}_n(\epsilon)\}|>\epsilon| C_r(y^0_{-K_n+1})|.
$$
Since $x_{-\infty}^{\infty}$ is typical, for large $N$, $|I_N({\tilde M}_n(\epsilon))|\le  \epsilon^2 N$ and 
\begin{eqnarray*}
\epsilon^2N&\ge& 
\sum_{y^0_{-K_n+1}\in{\cal X}^{K_n}}\sum_{r=1}^{R(y^0_{-K_n+1})} 
| \{ l\in C_r(y^0_{-K_n+1}): T^l x_{-\infty}^{\infty} \in {\tilde M}_n(\epsilon)\}|\\
&\ge&\epsilon 
\sum_{y^0_{-K_n+1}\in{\cal X}^{K_n}}\sum_{r=1}^{R(y^0_{-K_n+1})}  | C_r(y^0_{-K_n+1})|.
\end{eqnarray*}
Let $$I_N(B_n)=\{ n\le l\le N-1: T^{l-n}x_{-\infty}^{\infty}\in B_n\}.$$
But those $n\le l\le N-1$, such that $T^{l-n}x_{-\infty}^{\infty}\in B_n$ are  covered by this union - thus
$$
 \epsilon |I_N(B_n)|\le \epsilon^2 N
$$ 
and thus 
$$P(B_n)= \lim_{N\to\infty} {|I_N(B_n)|\over N}\le \epsilon
$$
since $x_{-\infty}^{\infty} $ was typical. 
The proof of the Theorem is complete.

\section{Appendix}

\begin{lemma} \label{applemma}
Let $\{X_n\}$ be stationary and ergodic.
For an arbitrary set $A$ measurable with respect to $\sigma(X^n_0)$,
the probability of the event 
$$
{\tilde A}(n,D)=\{x^{\infty}_{-\infty}\in A : \sum_{i=0}^{n-1} 1_{A}(T^{i} x^{\infty}_{-\infty})<D\}
$$
is not greater than $D/n$. 
\end{lemma} 
{\sc Proof} 
Fix a typical orbit $x^{\infty}_{-\infty}$. 
Let 
$$
I_N({\tilde A}(n,D))=\{ n\le l\le N-1: T^{l}x_{-\infty}^{\infty}\in {\tilde A}(n,D)\}.
$$
We make a disjoint cover.
Let $\dots,l_2<l_1$ denote the elements of $I_N({\tilde A}(n,D))$.
 Set 
$E_r=\emptyset$ and for $r=1,2,\dots$, define $E_r$ recursively. 
Let $l$ denote the largest element of $I_N({\tilde A}(n,D))$ such that 
$l\not \in \bigcup_{{r\prime}<r} E_{{r\prime}}$ if there is such and   let
$$E_r= \{ l-n\le l_i\le l:\ \mbox{for $i=1,2,\dots$.}\}
$$
Now let $R$ denote the largest $r$ for which $E_r$ has been defined. Since the cover is disjoint, 
$R (n+1)\le N$. Then clearly, 
$$
 {I_N({\tilde A}(n,D))\over N}\le  {R D\over R (n+1)}\le {D\over (n+1)}
$$
and the left hand side tends to $P({\tilde A}(n,D))$.
The proof of Lemma \ref{applemma} is complete. 


\end{document}